\documentclass[11pt,a4paper]{amsart}
\usepackage[T1]{fontenc}
\usepackage{lmodern,microtype}
\usepackage[margin=28mm]{geometry}
\usepackage{amsmath,amssymb,amsthm,mathtools,booktabs,array,longtable,enumerate}
\usepackage[hidelinks]{hyperref}
\hypersetup{pdftitle={Normal edge-transitive bi-Cayley graphs of prime valency over nonabelian simple groups}}
\newtheorem{theorem}{Theorem}[section]
\newtheorem{lemma}[theorem]{Lemma}

\newtheorem{corollary}[theorem]{Corollary}
\theoremstyle{definition}

\numberwithin{equation}{section}
\DeclareMathOperator{\Aut}{Aut}
\DeclareMathOperator{\Out}{Out}
\DeclareMathOperator{\Inn}{Inn}

\DeclareMathOperator{\BCay}{BCay}
\DeclareMathOperator{\Cay}{Cay}

\DeclareMathOperator{\Sp}{Sp}
\DeclareMathOperator{\supp}{supp}

\newcommand{\Fp}{\mathbb F_p}
\newcommand{\Alt}{\mathrm A}

\newcommand{\normal}{\trianglelefteq}

\title[Prime-valent bi-Cayley graphs]{Normal edge-transitive bi-Cayley graphs of prime valency  over nonabelian simple groups}
\thanks{Supported by the National Natural Science Foundation of China (12471328, 12331013) and the Fundamental Research Funds for the Central Universities.}

\author[Lu]{Zai Ping Lu}
\address{
Z.P. Lu\\
Center for Combinatorics, LPMC,
Nankai University,
Tianjin 300071, China
}
\email{lu@nankai.edu.cn}
\author[Wang]{Hai Ze Wang}
\address{
H.Z. Wang\\
Center for Combinatorics, LPMC,
Nankai University,
Tianjin 300071, China
}
\email{wanghaize0218@163.com}

\begin{document}

\begin{abstract}
Let $\Gamma$ be a connected graph of prime valency admitting a nonabelian
simple group $R$ that acts semiregularly on the vertices with exactly two
orbits. We prove that if the normalizer of $R$ in $\Aut(\Gamma)$ is
edge-transitive, then $\Gamma$ is bipartite and $R$ is characteristic in
$\Aut(\Gamma)$.  We derive consequences for Cayley
graphs and a criterion for semisymmetry. As an application, we construct
bi-Cayley graphs over $\Alt_p$ from the conjugates of a $3$-cycle under a
cyclic group of order $p$, where $p>5$ is prime. These graphs are
semisymmetric precisely when the support of the $3$-cycle is not a
three-term arithmetic progression in $\Fp$.

\vskip 5pt

\end{abstract}

\date{}

\maketitle

\parskip 5pt

\section{Introduction}\label{sec:introduction}

Throughout this paper all groups are finite. The graphs under consideration
are finite, simple and undirected, except for the auxiliary Cayley digraphs
introduced in Section~\ref{sec:consequences}, where loops are allowed.
A graph $\Gamma$ is a \emph{bi-Cayley graph} over a group $R$ if
$\Aut(\Gamma)$ contains a subgroup isomorphic to $R$ that acts
semiregularly on the vertices with exactly two orbits. We identify $R$
with this subgroup. The graph is \emph{normal}, with respect to $R$,
if $R\normal\Aut(\Gamma)$.

Two related conventions for bi-Cayley graphs occur in the literature.
In the bipartite formulation used by Du and Xu \cite{DuXu}, the defining
group acts regularly on each of the two parts, and these parts are
independent sets. The broader definition used in \cite{CZFF,PZsym,ZF}
allows edges within either of the two vertex-orbits. We adopt the broader
definition above, and reserve the notation $\BCay(R,S)$ for the bipartite
construction described in Section~\ref{sec:consequences}.

Praeger \cite{Praeger99} introduced normal edge-transitivity for Cayley
graphs. For $\Gamma=\Cay(R,S)$, this means that the normalizer of the
right multiplication subgroup $\widehat R$ is edge-transitive. This
normalizer is $\widehat R\Aut(R,S)$, where
$\Aut(R,S)=\{\alpha\in\Aut(R)\mid S^\alpha=S\}$.
Among the five questions posed in \cite{Praeger99}, Question~4 asks
what $\Aut(R,S)$ reveals about the full automorphism group and how large
the latter can be relative to this normalizer. Praeger proved that every
connected undirected cubic normal edge-transitive Cayley graph over a
nonabelian simple group is normal; in this case,
$\Aut(\Gamma)=\widehat R\Aut(R,S)$, see
\cite[Theorem~1]{Praeger99}. Zhang, Feng, Yin and Zhou subsequently
extended this normality result to connected normal edge-transitive
Cayley graphs of arbitrary prime valency over nonabelian simple groups
\cite{ZFYZ}.

For a bi-Cayley graph over $R$, we similarly use the term
\emph{normal edge-transitive} when $N_{\Aut(\Gamma)}(R)$ is
edge-transitive. This condition is called almost-normal edge-transitivity
in \cite{CZFF}. Lu \cite{Lu03} determined the structure of the
normalizer in the bipartite setting, and Zhou and Feng \cite{ZF}
determined it for bi-Cayley graphs in the broader sense. These descriptions
make the corresponding version of Praeger's question natural for
bi-Cayley graphs: when does the normalizer determine the full automorphism
group? Normal edge-transitivity alone does not require $R$ to be normal
in $\Aut(\Gamma)$.

Bi-Cayley graphs over nonabelian simple groups have been studied under
several symmetry assumptions. Pan and Zhang investigated cubic symmetric
graphs in \cite{PZsym} and cubic semisymmetric graphs in \cite{PZsemi}.
In particular, they proved that a connected cubic semisymmetric bi-Cayley
graph over a nonabelian simple group is normal whenever the normalizer
of that group is edge-transitive \cite[Corollary~1.2]{PZsemi}.
Prime-valent symmetric bi-Cayley graphs were studied under a solvability
assumption on vertex stabilizers in \cite{ZPWH}, and under
bi-quasiprimitivity assumptions in \cite{JLL}.
Here we assume edge-transitivity of the normalizer and prime valency,
and prove that the defining simple group is characteristic in the full
automorphism group. Our main result is the following.

\begin{theorem}\label{thm:main}
Let $\Gamma$ be a connected graph of prime valency, and let
$R\leqslant\Aut(\Gamma)$. Suppose that $R$ is a nonabelian simple group
and acts semiregularly with exactly two orbits on the vertex set. If
$N_{\Aut(\Gamma)}(R)$ acts transitively on the edge set of $\Gamma$, then
$\Gamma$ is bipartite, and $R$ is a characteristic subgroup of
$\Aut(\Gamma)$.
\end{theorem}

The theorem does not require vertex-transitivity or a solvability
assumption on vertex stabilizers. Its proof combines local actions of
prime degree with an analysis of minimal normal subgroups. The regular
action of $R$ on each part gives exact factorizations by $R$ and vertex
stabilizers, and the classifications in \cite{BL,LWX} exclude the
relevant simple overgroups. The characteristic-subgroup conclusion also
shows that $R$ is preserved by every abstract automorphism of
$\Aut(\Gamma)$.

Our application concerns semisymmetric graphs, that is, regular graphs
which are edge-transitive but not vertex-transitive. Their study goes
back to Folkman's constructions and problems \cite{Folkman}.
Subsequent work includes both constructions of infinite families and
classification results under restrictions on the order, valency or group
action. Du and Xu classified semisymmetric graphs of order $2pq$, where
$p$ and $q$ are distinct primes \cite{DuXu}. Conder, Malni\v{c},
Maru\v{s}i\v{c} and Poto\v{c}nik determined all connected cubic
semisymmetric graphs on at most $768$ vertices \cite{768-ss}, and Li
and Lu classified the connected pentavalent semisymmetric graphs of
square-free order \cite{LiLu}. In a setting closer to the present
construction, Lu, Wang and Xu constructed infinite families of cubic
semisymmetric bi-Cayley graphs over alternating groups \cite{LuWangXu06}.

We apply Theorem~\ref{thm:main} to $\BCay(R,S)$, where $R$ is a
nonabelian simple group and $S\subseteq R$. This graph has vertex set
$R\times\{+,-\}$, with $(x,+)$ adjacent to $(y,-)$ if and only if
$yx^{-1}\in S$. It is connected if and only if
$R=\langle S^{-1}S\rangle$. If this condition holds, $|S|$ is prime
and $\Aut(R,S)$ is transitive on $S$, then the full automorphism group
is the normalizer of $\widehat R$. Semisymmetry is then equivalent to
the absence of $x\in R$ and $\alpha\in\Aut(R)$ satisfying
$S^{-1}=xS^\alpha$.

For an explicit application, let $p>5$ be a prime, label the points of
the natural action of $\Alt_p$ by $\Fp$, and let $\rho$ be the translation
$u\mapsto u+1$. Given a $3$-cycle $\delta\in\Alt_p$, put
\[
\Delta_3=\{\delta^{\rho^i}\mid0\leqslant i\leqslant p-1\}.
\]
We prove that $\BCay(\Alt_p,\Delta_3)$ is connected and has valency $p$,
and that it is semisymmetric if and only if $\supp(\delta)$ is not a
three-term arithmetic progression in $\Fp$. This gives a criterion
in terms of the three points moved by $\delta$. In particular, choosing
$\delta=(0\,1\,3)$ gives a semisymmetric graph with $2|\Alt_p|=p!$ vertices for every
prime $p>5$.

The paper is organized as follows. Section~\ref{sec:local} records facts
about locally arc-transitive graphs. Section~\ref{sec:group} establishes
the group-theoretic lemmas used in the proof of Theorem~\ref{thm:main},
which is given in Section~\ref{sec:proof}. Section~\ref{sec:consequences}
derives consequences for Cayley and bi-Cayley graphs.
Section~\ref{sec:examples} gives the construction of semisymmetric graphs
over alternating groups.

\vskip 20pt

\section{Locally transitive graphs}\label{sec:local}

Let $\Gamma$ be a graph, and denote by $V(\Gamma)$, $E(\Gamma)$ and
$\Aut(\Gamma)$ its vertex set, edge set and automorphism group, respectively.
A subgroup $X\leqslant\Aut(\Gamma)$ is called a \emph{vertex-transitive}
or \emph{edge-transitive} group of $\Gamma$ if it acts transitively on
$V(\Gamma)$ or $E(\Gamma)$, respectively. For $v\in V(\Gamma)$, write
$X_v$ for its stabilizer in $X$, and $X_v^{\Gamma(v)}$ for the induced
permutation group on the neighbourhood $\Gamma(v)$. We call $X$
\emph{locally arc-transitive} or \emph{locally primitive} on $\Gamma$ if
$X_v^{\Gamma(v)}$ is transitive or primitive, respectively, for every
$v\in V(\Gamma)$. The group $X$ is \emph{arc-transitive} if it is both
vertex-transitive and locally arc-transitive. Note that a locally arc-transitive
group of a connected graph is edge-transitive.

If $\Gamma$ is connected and bipartite, let $X^+$ denote the subgroup of
$X$ that fixes each part of the bipartition setwise. Then
$|X:X^+|\leqslant2$. The following lemma collects some facts about
locally arc-transitive graphs.

\begin{lemma}\label{lem:edge-tran}
Let $\Gamma$ be a connected graph, let $X\leqslant\Aut(\Gamma)$, and let
$u\in V(\Gamma)$ and $w\in\Gamma(u)$.
\begin{enumerate}[\rm (1)]
\item If $X$ is edge-transitive but not vertex-transitive on $\Gamma$,
then $X$ is locally arc-transitive on $\Gamma$, and $\Gamma$ is bipartite
with its two parts equal to the $X$-orbits on $V(\Gamma)$.
\item If $X_u^{\Gamma(u)}$ and $X_w^{\Gamma(w)}$ are both transitive,
then $\langle X_u,X_w\rangle$ is locally arc-transitive on $\Gamma$.
If, in addition, $\Gamma$ is bipartite, then
$X^+=\langle X_u,X_w\rangle$.
\item If $\Gamma$ has prime valency $p$ and $X$ is locally
arc-transitive on $\Gamma$, then $p$ is the largest prime divisor of
$|X_u|$, and $p$ does not divide $|X_{uw}|$.
\end{enumerate}
\end{lemma}
\begin{proof}
Suppose first that $X$ is edge-transitive but not vertex-transitive.
Put $U=u^X$ and $W=w^X$. Edge-transitivity and connectedness give
$V(\Gamma)=U\cup W$, and the assumption that $X$ is not
vertex-transitive gives $U\cap W=\emptyset$. For every edge
$\{v_1,v_2\}$, some $x\in X$ satisfies
$\{u,w\}^x=\{v_1,v_2\}$. Thus one end of each edge lies in $U$ and
the other lies in $W$, so $(U,W)$ is the bipartition of $\Gamma$.
Let $v\in V(\Gamma)$ and $w_1,w_2\in\Gamma(v)$. There exists $y\in X$
such that $\{v,w_1\}^y=\{v,w_2\}$. Since $X$ preserves each part,
$y$ fixes $v$ and maps $w_1$ to $w_2$. Hence $X_v$ is transitive on
$\Gamma(v)$, and (1) of the lemma follows.

Part (2) is recorded in \cite[Exercise~3.8]{Weiss} without a proof.
Put $Y=\langle X_u,X_w\rangle$. For any vertex $v$, connectedness
gives a walk $v_0=u,v_1=w,\ldots,v_\ell=v$. We prove inductively
that $X_{v_i}\leqslant Y$ and that $X_{v_i}^{\Gamma(v_i)}$ is
transitive. The assertions hold for $i=0,1$. Suppose they hold for
$i-1$ and $i$. Choose $x\in X_{v_i}$ such that
$v_{i-1}^x=v_{i+1}$. Then
$X_{v_{i+1}}=X_{v_{i-1}}^x\leqslant Y$, and its local action is
transitive by conjugacy. Consequently, $Y_v=X_v$ and
$Y_v^{\Gamma(v)}$ is transitive for every $v\in V(\Gamma)$.
Thus $Y$ is locally arc-transitive. If $\Gamma$ is bipartite, then
$Y\leqslant X^+$, and both groups are edge-transitive. Their
stabilizers of $\{u,w\}$ are both $X_{uw}=X_u\cap X_w$, so
\[
|Y:X_{uw}|=|E(\Gamma)|=|X^+:X_{uw}|.
\]
It follows that $Y=X^+$, proving (2).

For (3), let $K$ be the kernel of the action of $X_u$ on $\Gamma(u)$.
Suppose that a prime $q\geqslant p$ divides $|K|$, and let $Q$ be a
 Sylow $q$-subgroup of $K$. Then $Q$ fixes $u$ and every
vertex in $\Gamma(u)$. Whenever $Q$ fixes adjacent vertices $v,w$,
it acts on the invariant set $\Gamma(w)\setminus\{v\}$, which has
$p-1<q$ elements. Every orbit has $q$-power size, so $Q$ fixes this
set pointwise. By connectedness, $Q$ fixes every vertex of $\Gamma$,
a contradiction. Thus every prime divisor of $|K|$ is less than $p$.
Since $X_u^{\Gamma(u)}$ is a transitive subgroup of ${\rm S}_p$,
its order is divisible by $p$ but not by $p^2$, and has no prime
divisor greater than $p$. The same assertions hold for $|X_u|$.
Finally, $|X_u:X_{uw}|=p$, and hence $p$ does not divide $|X_{uw}|$.
\end{proof}

\begin{lemma}\label{lem:prime}
Let $\Gamma$ be a connected graph of prime valency $p$, let
$\{u,w\}\in E(\Gamma)$, and let $X\leqslant\Aut(\Gamma)$. If $p$
divides both $|X_u|$ and $|X_w|$, then $X$ is locally arc-transitive
on $\Gamma$.
\end{lemma}
\begin{proof}
By Lemma~\ref{lem:edge-tran}(2), it suffices to show that
$X_u^{\Gamma(u)}$ and $X_w^{\Gamma(w)}$ are both transitive.
Suppose that $X_u^{\Gamma(u)}$ is not transitive, and let $P$ be a
Sylow $p$-subgroup of $X_u$. Then every $P$-orbit on $\Gamma(u)$ has
$p$-power size. An orbit of size $p$ would make
$X_u^{\Gamma(u)}$ transitive, so $P$ fixes $\Gamma(u)$ pointwise.
For each $v\in\Gamma(u)$, the group $P$ fixes $u$ and $v$ and
therefore acts on the $(p-1)$-set $\Gamma(v)\setminus\{u\}$.
All orbits on this set have size one. Continuing along paths in
$\Gamma$, we find that $P$ fixes every vertex. Thus $P=1$, contrary
to $p\mid|X_u|$. The same argument applies at $w$, and the lemma
follows.
\end{proof}

\begin{lemma}\label{lem:bipartition}
Let $\Gamma$ be a connected graph, and let
$R\leqslant A=\Aut(\Gamma)$. Suppose that $R$ has exactly two
orbits, say $U$ and $W$, on $V(\Gamma)$, and that $N_A(R)$ is
edge-transitive on $\Gamma$. Then $\Gamma$ is bipartite with parts
$U$ and $W$. If $N_{A^+}(R)$ is not locally arc-transitive on
$\Gamma$, then $\Gamma$ is either not regular or has even valency.
\end{lemma}
\begin{proof}
Note that the normalizer $N_A(R)$ permutes the two $R$-orbits. Suppose that
$\Gamma$ has an edge $\{v,v'\}$ with both ends in the same orbit,
say $U$. Since $\Gamma$ is connected, it also has an edge
$\{u,w\}$ with $u\in U$ and $w\in W$. Choose $g\in N_A(R)$ such
that $\{v,v'\}^g=\{u,w\}$. Both $u$ and $w$ then belong to the
same $R$-orbit $U^g$, a contradiction. Hence no edge has both ends
in the same orbit, and so $(U,W)$ is the bipartition of $\Gamma$. In particular,
$R\leqslant A^+$. Put $N=N_A(R)$. Then $N^+=N_{A^+}(R)$.

Suppose that $N^+$ is not locally arc-transitive. By
Lemma~\ref{lem:edge-tran}(1), it is not edge-transitive, so
$|N:N^+|=2$. Hence $N^+$ has exactly two orbits $E_1,E_2$ on
$E(\Gamma)$, and some $g\in N\setminus N^+$ interchanges them.
Since $R\leqslant N^+$ is transitive on both of $U$ and $W$,
each $E_i$ is $R$-invariant and defines a spanning subgraph
$\Gamma_i$ of $\Gamma$. Suppose now that $\Gamma$ is regular.
Counting edges between the two parts gives $|U|=|W|$. As $R$ is
transitive on each part, each $\Gamma_i$ has constant valency on
each part; counting its edges shows that these two valencies are
equal. Thus $\Gamma_1$ and $\Gamma_2$ are regular. The element $g$
gives an isomorphism between them, so they have the same valency.
Since $\Gamma$ is the edge-disjoint union of these two subgraphs, $\Gamma$   has even valency.
\end{proof}

\begin{lemma}\label{lem:normal-local}
Let $\Gamma$ be a connected bipartite graph, and let
$Y\normal X\leqslant A^+$, where $A=\Aut(\Gamma)$. Suppose that
$X$ is locally primitive on $\Gamma$, and that $Y$ is transitive
on both parts. Then $Y$ is either locally arc-transitive on
$\Gamma$ or regular on both parts.
\end{lemma}
\begin{proof}
For every $v\in V(\Gamma)$, we have $Y_v\unlhd X_v$. By local
primitivity, $Y_v^{\Gamma(v)}$ is either trivial or transitive.
If it is transitive for some $v$, then the transitivity of $Y$ on
the part containing $v$ implies that $Y$ is edge-transitive.
Lemma~\ref{lem:edge-tran}(1) then makes $Y$ locally arc-transitive.
Otherwise, $Y_v^{\Gamma(v)}=1$ for every $v\in V(\Gamma)$.
For adjacent vertices $u,w$, this gives $Y_u\leqslant Y_w$ and
$Y_w\leqslant Y_u$, so $Y_u=Y_w$. By connectedness, all vertices
have the same stabilizer in $Y$. This stabilizer fixes every
vertex and is therefore trivial. Thus $Y$ is semiregular on
$V(\Gamma)$ and regular on each part.
\end{proof}

Every transitive group of prime degree is primitive. The following
two lemmas consequently follow from \cite[Lemmas~5.1 and 5.2]{GLP}.

\begin{lemma}\label{lem:intransitive}
Let $\Gamma$ be a connected bipartite graph of prime valency, and
let $K\normal X\leqslant A^+$, where $A=\Aut(\Gamma)$. Suppose
that $K$ is intransitive on both parts. If $X$ is locally
arc-transitive on $\Gamma$, then $K$ is semiregular on $V(\Gamma)$.
\end{lemma}

\begin{lemma}\label{lem:faithful}
Let $\Gamma$ be a connected bipartite graph of prime valency $p$,
and let $X\leqslant A^+$, where $A=\Aut(\Gamma)$. If $X$ is
locally arc-transitive on $\Gamma$, then either
$\Gamma\cong K_{p,p}$ or $X$ acts faithfully on each part.
\end{lemma}

\vskip 20pt

\section{Several technical lemmas}\label{sec:group}

The following result is a consequence of \cite[Corollary]{Ka};
see also \cite[Lemma~A.2]{EGS}.

\begin{lemma}\label{eq:fixedpoint}
Let $L$ be a finite nonabelian simple group, and let
$1\neq\alpha\in\Aut(L)$. Then $|L:C_L(\alpha)|$ is not a prime
power.
\end{lemma}

\begin{lemma}\label{lem:diagonal}
Let $X$ be a locally primitive group of a connected bipartite
graph $\Gamma$ of prime power valency $p^e$. Suppose that
$N=L\times K\unlhd X$, where $L$ and $K$ are finite nonabelian
simple groups. Then $N\leqslant X^+$. Moreover, if $L$ is regular
on one part of $\Gamma$, then $K$ is not regular on the other part.
\end{lemma}
\begin{proof}
Since $L$ and $K$ are nonabelian simple, $N$ has no subgroup of
index two. Hence $N\leqslant X^+$.

Let $\{u,w\}\in E(\Gamma)$, and write $U=u^{X^+}$ and
$W=w^{X^+}$ for the two parts. Suppose that $L$ is regular on $U$
and $K$ is regular on $W$. Since $\Gamma$ is regular, its two
parts have the same size, so
\[
|L|=|U|=|W|=|K|.
\]
This order is not a prime power, and hence
$\Gamma\not\cong K_{p^e,p^e}$. By \cite[Lemma~5.2]{GLP}, $X^+$
acts faithfully on both parts. Note that the centralizer of a transitive
permutation group is semiregular; see
\cite[page~109, Theorem~4.2A]{Dixon}. Since $K$ centralizes $L$,
it is semiregular on $U$ and therefore regular on $U$. Similarly,
$L$ is regular on $W$.

For each $x\in L$, let $\hat{x}\in K$ be the unique element with
$u^x=u^{\hat{x}}$. Since $L$ and $K$ commute, for $x,y\in L$ we have
\[
u^{\hat{x}\hat{y}}=(u^x)^{\hat{y}}
=u^{x\hat{y}}=u^{\hat{y}x}=u^{yx}=u^{\widehat{yx}}.
\]
The regularity of $K$ gives
$\hat{x}\hat{y}=\widehat{yx}$ and
$\widehat{x^{-1}}=\hat{x}^{-1}$. Thus
\[
\phi:L\longrightarrow K,\,\,\,\, x\longmapsto\hat{x}^{-1},
\]
is an isomorphism. Noting that $x\hat{x}^{-1}$ fixes $u$, since
$N=LN_u$ and $L\cap N_u=1$, we have $|N_u|=|K|=|L|$, and hence
\[
N_u=\{x\hat{x}^{-1}\mid x\in L\}.
\]
Likewise, for each $x\in L$, choose $\tilde{x}\in K$ such that
$w^x=w^{\tilde{x}}$. Then
\[
\psi:L\longrightarrow K,\,\,\,\, x\longmapsto\tilde{x}^{-1},
\]
is an isomorphism and
\[
N_w=\{x\tilde{x}^{-1}\mid x\in L\}.
\]
The uniqueness of the two coordinates in $L\times K$ gives
\[
N_{uw}=\{x\tilde{x}^{-1}\mid x\in L,\ \phi(x)=\psi(x)\}.
\]
Define $\alpha\in\Aut(L)$ explicitly by
$\alpha(x)=\phi^{-1}(\psi(x))$. Then
$\phi(x)=\psi(x)$ if and only if $x\in C_L(\alpha)$, so
\[
|N_u:N_{uw}|=|L:C_L(\alpha)|.
\]
If $\alpha=1$, this index is one. If $\alpha\neq1$, it is not
a prime power by Lemma~\ref{eq:fixedpoint}. Thus in either case
it is not $p^e$. On the other hand, $N$ is transitive and not
regular on either part. Applying Lemma~\ref{lem:normal-local}
to $N\unlhd X^+$ shows that $N$ is locally arc-transitive, and
hence $|N_u:N_{uw}|=p^e$. This contradiction proves the lemma.
\end{proof}

For a positive integer $n$ and a prime $p$, let $n_p$ denote the
highest power of $p$ dividing $n$, and let $\pi(n)$ be the set
of prime divisors of $n$. For a finite group $X$, write
$\pi(X)=\pi(|X|)$.

\begin{lemma}\label{lem:transfer}
Let $X$ be a finite group, and let $p\in\pi(X)$. Suppose that $X$
has a nontrivial perfect subgroup $R$ and a normal subgroup $K$
of order $p$ such that $R\cap K=1$ and $|X:R|_p=p$. Then there
exists $L\unlhd X$ such that $K\cap L=1$ and $R\leqslant L$.
\end{lemma}
\begin{proof}
The conjugation of $R$ on $K$ induces an embedding of $R/C_R(K)$ into the abelian
group $\Aut(K)$. Since $R$ is perfect, $R=C_R(K)$. Put
$C=C_X(K)$. Then  $RK=R\times K\leqslant C \unlhd X$, and the
assumption $|X:R|_p=p$ gives
\[
|RK|_p=|C|_p=|X|_p.
\]
If $P$ is a Sylow $p$-subgroup of $R$, then $P\times K$ is a
Sylow $p$-subgroup of $C$. Thus $K$ has a complement in a Sylow
$p$-subgroup of $C$. By Gasch\"utz's theorem
\cite[page~31, (10.4)]{Aschbacher}, it has a complement $H$ in
$C$. Since $K\leqslant Z(C)$, we have $C=K\times H$.
Let $L=C'$. Then $L=H'$, $K\cap L=1$ and $L\unlhd X$.
Moreover, $R=R'\leqslant C'=L$, as required.
\end{proof}

\begin{lemma}\label{lem:factorization}
Let $T$ be a finite nonabelian simple group with an exact
factorization $T=RB$, where $R,B<T$ and $R\cap B=1$. Let
$p=\max\pi(B)$. If $R$ is nonabelian simple and $|B|_p=p$,
then $|N_T(R):R|_p=1$.
\end{lemma}
\begin{proof}
If $p=2$, then $|B|=2$, so $R$ has index two in $T$, a
contradiction. Thus $p$ is odd. Suppose that
$|N_T(R):R|_p\neq1$. Then $|N_T(R):R|\geqslant3$. Since $R$
is a proper subgroup of the simple group $T$, its normalizer
is proper in $T$, and hence $R$ is not maximal in $T$.
We apply the exact-factorization classifications in
\cite{BL,LWX}, retaining the entries in which $T$ and a
nonmaximal factor are nonabelian simple. In the notation of
those papers, put $G=T$ and $\{H,K\}=\{R,B\}$.

\emph{Case 1.} Suppose that one factor, say $H$, is solvable.
Then $B=H$ and $R=K$, and case (ii), (iv) or (v) of
\cite[Theorem~3]{BL} occurs. In case (ii), $T=G={\rm A}_n$ and
$R=K={\rm A}_{n-2}$ is the pointwise stabilizer of a two-element
set in the natural action. Thus $N_T(R)\cong{\rm S}_{n-2}$,
and $|N_T(R):R|=2$, a contradiction. In case (iv),
$T=G={\rm PSp}_{2m}(q)$ and $R=K=\Omega^-_{2m}(q)$, where $q$
is even and $m\geqslant3$. Here $R$ has index two in the maximal
subgroup ${\rm O}^-_{2m}(q)$ of $T$, so again
$|N_T(R):R|=2$, a contradiction.

In case (v), the factorization is listed in \cite[Table~4]{BL}.
Under the conditions that $G$ and $K$ are nonabelian simple and
$K$ is not maximal in $G$, the remaining possibility is
$T=G={\rm A}_n$ and $R=K={\rm A}_{n-2}$, with
$n\in\{25,49,121,529\}$. As above, $R$ is the pointwise
stabilizer of a two-element set, and $|N_T(R):R|=2$, a
contradiction.

\emph{Case 2.} Suppose that both factors are nonsolvable.
One of them is $R$, which is nonabelian simple and not maximal
in $G$. Row 15 of \cite[Table~1.1]{LWX} gives the possibility
\[
T=\Sp_6(2),\,\,\,\, R=\mathrm{SU}_3(3),
\,\,\,\, B={\rm S}_5\times C_2.
\]
Here $R$ is normal and of index two in the maximal subgroup
$\mathrm G_2(2)\cong\mathrm{SU}_3(3).2$ of $T$; see
\cite{Atlas}. Since $N_T(R)<T$, the maximality of $G_2(2)$ in $T$ gives
$N_T(R)=\mathrm G_2(2)$ and hence $|N_T(R):R|=2$.
This contradicts the assumption that the odd prime $p$ divides
$|N_T(R):R|$.

After excluding this case, the relevant cases in
\cite[Theorem~1.2 and Table~1.1]{LWX} come from Rows 1--4,
6, 7 and 9--11 of the table. More precisely, the remaining
possibilities for $T,R,B$ are as follows:
\begin{enumerate}[(i)]
\item $T={\rm A}_n$, $R$ is regular in the natural action of
${\rm A}_n$, and $B={\rm A}_{n-1}$ is a point stabilizer;
\item $T={\rm A}_{q+1}$, $R=\mathrm{PSL}_2(q)$, and
$B={\rm A}_{q-2}$, where $q>4$ is even;
\item $T={\rm A}_{q+1}$, $R={\rm A}_{q-2}$ is the pointwise
stabilizer of a three-element set, and
$B=\mathrm{PSL}_2(q).(2,q-1)$, with the embeddings described in
\cite[Example~4.1]{LWX};
\item $T={\rm A}_{q+1}$, $R=\mathrm{PSL}_2(q)$, and
$B={\rm A}_{q-2}.2$, where $q\equiv3\pmod4$;
\item $T={\rm A}_9$, $R=\mathrm{PSL}_2(8)$, and
$B=({\rm A}_5\times3).2$;
\item $T={\rm A}_{11}$, $R={\rm A}_7$ is the pointwise
stabilizer of a four-element set, and $B={\rm M}_{11}$;
\item $T={\rm A}_{12}$, $R={\rm A}_7$ is the pointwise
stabilizer of a five-element set, and $B={\rm M}_{12}$;
\item $T={\rm A}_n$, $R={\rm A}_{n-2}$ is the pointwise
stabilizer of a two-element set, and $B$ is sharply
two-transitive, where $n\in\{121,841,3481\}$.
\end{enumerate}
We rule out these possibilities in turn. Since $B$ is
nonsolvable, $p\geqslant5$.

Suppose first that (iii) or one of (vi)--(viii) holds. Then
$R$ is the pointwise stabilizer of a $t$-element set in
${\rm A}_n$, where $2\leqslant t\leqslant5$. Consequently,
\[
N_T(R)=({\rm A}_{n-t}\times{\rm A}_t).2,
\,\,\,\, |N_T(R):R|=t!.
\]
Since $p\geqslant5$ divides this index, we must have $t=p=5$.
This leaves only (vii), where $B={\rm M}_{12}$. But then
$p=\max\pi({\rm M}_{12})=11$, a contradiction.

Suppose that (ii), (iv) or (v) holds. The group $R$ is doubly
transitive in the natural action of $T={\rm A}_n$, and hence
$C_T(R)=1$; see \cite[page~109, Theorem~4.2A]{Dixon}.
Thus $N_T(R)/R$ embeds in $\Out(R)$. Write $q=r^f$, where $r$
is prime; in case (v), take $q=8$. Then
\[
|N_T(R):R|\ \text{divides}\ (2,q-1)f.
\]
Our assumption gives $p\mid(2,q-1)f$. As $p\geqslant5$, it
follows that $p\mid f$, and hence $p\leqslant f$ and $f\geqslant5$.
In particular, case (v) is excluded. By Bertrand's postulate,
there is a prime $s$ with $f<s<2f$. Since
$q-2\geqslant2^f-2>2f$, we have $s<q-2$, and hence $s$ divides
$|{\rm A}_{q-2}|$. In cases (ii) and (iv), this order divides
$|B|$. It follows that $p\leqslant f<s\leqslant\max\pi(B)=p$,
a contradiction.

Finally, suppose that (i) holds. Here $|R|=n$, and the normalizer
of the regular subgroup $R$ in ${\rm S}_n$ is its holomorph.
Thus
\[
N_{{\rm S}_n}(R)/R\cong\Aut(R).
\]
Since $p$ divides $|N_{{\rm A}_n}(R):R|$, it divides
$|\Aut(R)|$, so $R$ has an automorphism $\alpha$ of order $p$.
Note that the action of $\langle\alpha\rangle$ on $R$ has a fixed point
and a nontrivial orbit of length $p$. Since
$|{\rm A}_{n-1}|_p=p$ and $p$ is odd, we have
$p\leqslant n-1<2p$. Hence there is exactly one orbit of length
$p$ and exactly $n-p$ fixed points. Therefore,
$|C_R(\alpha)|=n-p$ divides $n$. If $n-p>1$, any prime divisor
of $n-p$ divides both $n$ and $p$, and hence equals $p$.
Since $n\leqslant2p$, this forces $n=2p$, contrary to the
nonabelian simplicity of $R$. Thus $n=p+1$. In particular,
$|R|$ is even and $R$ contains an involution. As
$\langle\alpha\rangle$ is transitive on $R\setminus\{1\}$,
every nonidentity element of $R$ is an involution. This makes
$R$ abelian, a final contradiction.
\end{proof}

\vskip 20pt

\section{Proof of Theorem~\ref{thm:main}}\label{sec:proof}

Let $\Gamma$ and $R$ satisfy the hypotheses of
Theorem~\ref{thm:main}. Put $A=\Aut(\Gamma)$ and $N=N_A(R)$,
and let $p$ be the valency of $\Gamma$. A connected graph of
valency two is a cycle and has a solvable automorphism group.
Thus $p$ is odd. By Lemma~\ref{lem:bipartition}, $\Gamma$ is
bipartite and $N^+=N_{A^+}(R)$ is locally arc-transitive, and
the two $R$-orbits, denoted by $U$ and $W$, form the bipartition.
Since $R$ is regular on each part, $|U|=|W|=|R|$. The valency
gives $|U|\geqslant p$, and equality is impossible because $R$
is nonabelian simple. Hence $|R|>p$ and
$\Gamma\not\cong K_{p,p}$. By Lemma~\ref{lem:faithful}, $A^+$
acts faithfully on each part. By Lemma~\ref{lem:edge-tran}(1),
a subgroup of $A^+$ is edge-transitive if and only if it is
locally arc-transitive.

\begin{lemma}\label{lem:complement}
Suppose that $A^+$ has a simple subgroup $T$ such that
$R\leqslant N_{A^+}(T)$ and $RT$ is edge-transitive on $\Gamma$.
Then either $R<T$, or $|T|_p=p$, $RT=R\times T$ and $|T|<|R|$.
\end{lemma}
\begin{proof}
Since $R\cap T\unlhd R$ and $R$ is simple, $R\cap T=1$ or $R$.
The edge-transitivity of $RT$ gives $|RT|\geqslant p|U|=p|R|$, so
$|T:R\cap T|\geqslant p$. Thus, if $R\cap T=R$ then $R<T$.

Suppose that $R\cap T=1$, and put $B=RT$. Since $R$ is regular
on $U$, we have $|B_u|=|T|$ for $u\in U$. By
Lemma~\ref{lem:edge-tran}(3), $|T|_p=|B_u|_p=p$.
Note that the conjugation of $B$ on $T$ induces a homomorphism
$\phi:B\longrightarrow\Aut(T)$ with kernel $C=C_B(T)$.
If $T$ is abelian, its simplicity and $|T|_p=p$ give $|T|=p$.
Then $\Aut(T)$ is abelian, so $\phi(R)=1$. Thus
$RT=R\times T$ and $|T|=p<|R|$, as required.

Suppose now that $T$ is nonabelian. Since $\phi(T)=\Inn(T)$,
\[
B/CT\cong\phi(B)/\phi(T)\leqslant\Out(T).
\]
The group $B/CT$ is a quotient of $B/T\cong R$, whereas
$\Out(T)$ is solvable. Hence $B=CT$. Since $C\cap T=Z(T)=1$,
we have $B=C\times T$ and $C\cong R$.

Suppose that $T$ is transitive on one part, say $U$. As $B$
acts faithfully on $U$, its subgroup $C$, which centralizes
$T$, is semiregular there; see
\cite[page~109, Theorem~4.2A]{Dixon}. Since $|C|=|R|=|U|$,
the group $C$ is regular on $U$. The same centralizer argument
then shows that $T$ is regular on $U$. By
Lemma~\ref{lem:diagonal}, neither $T$ nor $C$ is regular on
$W$. Neither can be transitive on $W$, since the same argument
would make both regular there. By \cite[Lemma~5.5]{GLP}, each of $T$ and $C$
has exactly $p$ orbits on $W$. These orbits have equal size
because $T,C\unlhd B$ and $B$ is transitive on $W$. For
$w\in W$, we therefore have
\[
|T:T_w|=\frac{|W|}{p}=\frac{|T|}{p},
\,\,\,\, |T_w|=p.
\]
Similarly, $|C_w|=p$. Since $T_w\times C_w\leqslant B_w$,
this contradicts $|B_w|_p=p$. Hence $T$ is intransitive on
both parts. Lemma~\ref{lem:intransitive} makes $T$ semiregular,
and consequently $|T|<|U|=|R|$.

Since $C\unlhd B$, either $R\leqslant C$ or $R\cap C=1$.
In the latter case, $\phi(R)\cong R$. Its image in the solvable
group $\Out(T)$ is trivial, so
$\phi(R)\leqslant\Inn(T)\cong T$. This contradicts
$|T|<|R|$. Therefore $R\leqslant C$, and $RT=R\times T$.
\end{proof}

\begin{proof}[Proof of Theorem~\ref{thm:main}]
We retain the notation introduced above. We first prove that
$R\unlhd A^+$. Suppose, to the contrary, that
\begin{itemize}
\item[($\ast$)] $R$ is not normal in $A^+$, and $T$ is a
minimal normal subgroup of $A^+$.
\end{itemize}
Since $R\cap T\unlhd R$, either $R\cap T=1$ or $R<T$.

\smallskip
\noindent\textbf{Claim 1.}
\emph{Under assumption ($\ast$), the group $T$ is simple and
$RT$ is locally arc-transitive on $\Gamma$.}

We have $R<RT\normal N^+T$. The group $N^+T$ is locally
primitive, since its local actions contain transitive groups
of prime degree. The subgroup $RT$ is transitive on both
parts and is not regular there, because it properly contains
$R$. Lemma~\ref{lem:normal-local} therefore makes $RT$ locally
arc-transitive. If $R\cap T=1$, then
$|(RT)_u|=|T|$ and $|T|_p=p$. As $T$ is characteristically
simple, it is a direct product of isomorphic simple groups.
The equality $|T|_p=p$ forces there to be only one factor.

Suppose that $R<T$. Write $T=T_1\times\cdots\times T_m$,
where the $T_i$ are isomorphic nonabelian simple groups.
Local transitivity gives $|T_u|_p=p$, so $p$ divides the order
of every $T_i$. Since $R\leqslant T$ normalizes each factor,
$R\cap T_i=1$ or $R$. If $m>1$, then $R\cap T_j=1$ for some
$j$. For every vertex $v$, the regularity of $R$ gives
$|(RT_j)_v|=|T_j|$, which is divisible by $p$. By
Lemma~\ref{lem:prime}, $RT_j$ is locally arc-transitive, and
Lemma~\ref{lem:complement} gives $|T_j|<|R|$.
Noting that at least one projection of $R$ onto a simple direct factor
$T_i$ is nontrivial and hence injective, we have
$|R|\leqslant|T_i|=|T_j|$, a contradiction. Consequently,
$m=1$, and the claim follows.

\smallskip
\noindent\textbf{Claim 2.}
\emph{Under assumption ($\ast$), the group $T$ is the unique
minimal normal subgroup of $A^+$.}

Suppose that $S$ is another minimal normal subgroup of $A^+$.
Then $S\cap T=1$ and $ST=S\times T$. By Claim 1, $S$ is
simple and both $RS$ and $RT$ are locally arc-transitive.
If $R\cap S=R\cap T=1$, Lemma~\ref{lem:complement} gives
\[
RS=R\times S,\,\,\,\, RT=R\times T,\,\,\,\, |S|_p=|T|_p=p.
\]
Thus $ST$ centralizes $R$, so $R\cap ST\leqslant Z(R)=1$.
For every vertex $v$, we have
$|(RST)_v|=|S||T|$ and hence $|(RST)_v|_p=p^2$.
Lemma~\ref{lem:prime} makes $RST$ locally arc-transitive,
contrary to Lemma~\ref{lem:edge-tran}(3).

We may therefore suppose that $R\cap S\neq1$. Then $R<S$,
and $R\cap T=1$ because $S\cap T=1$. Both $S$ and $ST$ are
normal in $A^+$, transitive on both parts, and not regular.
By Lemma~\ref{lem:normal-local}, they are locally
arc-transitive. For every vertex $v$, we have
$S=RS_v$, $R\cap S_v=1$ and
\[
|(ST)_v|=|S_v||T|.
\]
Also $|(RT)_v|=|T|$. Since $S$, $RT$ and $ST$ are locally
arc-transitive, Lemma~\ref{lem:edge-tran}(3) gives
$|S_v|_p=|T|_p=|(ST)_v|_p=p$, contradicting the displayed
equality. This proves Claim 2.

\smallskip
\noindent\textbf{Claim 3.}
\emph{Under assumption ($\ast$), we have $R<T$.}

Otherwise, Claim 1 and Lemma~\ref{lem:complement} give
$|T|_p=p$ and $RT=R\times T$. If $T$ is nonabelian, then
$T\cap C_{A^+}(T)=1$ and
$1\neq R\leqslant C_{A^+}(T)\unlhd A^+$, contrary to
Claim 2. Thus $T$ is abelian and $|T|=p$. For $u\in U$, the
regularity of $R$ gives $A^+=RA_u$ and $R\cap A_u=1$.
Consequently,
\[
|A^+:R|_p=|A_u|_p=p.
\]
By Lemma~\ref{lem:transfer}, there exists $L\unlhd A^+$ such
that $R\leqslant L$ and $T\cap L=1$, again contradicting
Claim 2. Hence $R<T$.

\smallskip
We now derive a contradiction from ($\ast$). Fix $u\in U$.
Since $R$ is regular on $U$, we have $T=RT_u$ and
$R\cap T_u=1$. Lemma~\ref{lem:edge-tran}(3) gives
\[
|A_u|_p=|N_u|_p=|T_u|_p=p,
\,\,\,\, p=\max\pi(T_u).
\]
Since $T\unlhd A^+$, we have $T_u\unlhd A_u$ and hence
$N_uT_u\leqslant A_u$. It follows that
\[
p^2=|N_u|_p|T_u|_p
=|N_uT_u|_p|N_u\cap T_u|_p,
\]
so $|N_u\cap T_u|_p=p$. Moreover,
$N_u\cap T_u\leqslant N_T(R)$ and
$R\cap(N_u\cap T_u)=1$. Thus $p$ divides $|N_T(R):R|$,
contrary to Lemma~\ref{lem:factorization} applied to the exact
factorization $T=RT_u$. Therefore $R\unlhd A^+$.

For $a\in A$, the subgroup $R^a$ is also a nonabelian simple
normal subgroup of $A^+$. Hence either $R^a=R$ or
$RR^a=R\times R^a$. In the latter case, both factors are
regular on both parts, since conjugation by $a$ preserves
this property. Applying Lemma~\ref{lem:diagonal} to
$R\times R^a\unlhd A^+$ gives a contradiction. Thus
$R^a=R$ for every $a\in A$, and $R\unlhd A$.

Finally, let $\alpha\in\Aut(A)$. Then $R^\alpha\unlhd A$.
Since $R^\alpha$ is nonabelian simple, it has no subgroup
of index two, and therefore $R^\alpha\leqslant A^+$.
If $R^\alpha\neq R$, then these two simple normal subgroups
of $A$ have trivial intersection and centralize one another.
As $A^+$ is faithful on each part and $R$ is regular there,
$R^\alpha$ is semiregular on each part. Since
$|R^\alpha|=|R|=|U|=|W|$, $R^\alpha$ is regular on both parts.
Lemma~\ref{lem:diagonal}, applied to
$R\times R^\alpha\unlhd A^+$, yields another contradiction.
Hence $R^\alpha=R$ for every $\alpha\in\Aut(A)$, and $R$ is
characteristic in $A$.
\end{proof}

\vskip 20pt

\section{Consequences}\label{sec:consequences}

Let $R$ be a finite group and $S\subseteq R$, possibly containing the
identity element $1$. The Cayley digraph $\Cay(R,S)$, with loops allowed,
has vertex set $R$ and arc set $\{(x,y)\mid yx^{-1}\in S\}$. If
$S=S^{-1}$, we regard $\Cay(R,S)$ as an undirected graph with edge set
$\{\{x,y\}\mid yx^{-1}\in S\}$, where $\{x,x\}$ is a loop when
$1\in S$. It is well known that $\Cay(R,S)$ is connected if and only if
$R=\langle S\rangle$. Moreover, $\Aut(\Cay(R,S))$ contains a
vertex-regular subgroup $\hat R$ induced by the right multiplication
action of $R$ on itself. Clearly, $\hat R\cong R$.

The associated bipartite bi-Cayley graph $\BCay(R,S)$ has vertex set
$R\times\{+,-\}$ and edge set
\[
\bigl\{\{(x,+),(y,-)\}\mid yx^{-1}\in S\bigr\}.
\]
When $S=S^{-1}$, this is the standard double cover of $\Cay(R,S)$.
By \cite[Lemma~2.3(iii)]{DuXu}, $\BCay(R,S)$ is connected if and only if
$S\ne\emptyset$ and $R=\langle S^{-1}S\rangle$.

We regard $\Aut(\Cay(R,S))$ as a subgroup of $\Aut(\BCay(R,S))$ via
\[
(x,\epsilon)^a=(x^a,\epsilon),\,\,\,\,
x\in R,\quad \epsilon\in\{+,-\},\quad a\in\Aut(\Cay(R,S)).
\]
If $X\leqslant\Aut(\Cay(R,S))$ is arc-transitive on $\Cay(R,S)$,
then $X$ is locally arc-transitive on $\BCay(R,S)$. Under this
identification, $\hat R$ is semiregular on $R\times\{+,-\}$ with
orbits $R\times\{+\}$ and $R\times\{-\}$, where
\[
(x,\epsilon)^{\hat h}=(xh,\epsilon),\,\,\,\,
x,h\in R,\quad \epsilon\in\{+,-\}.
\]

Conversely, let $\Gamma$ be a bipartite graph with bipartition $(U,W)$,
and suppose that $\Aut(\Gamma)$ contains a subgroup $R$ that is regular
on both $U$ and $W$. Choose $u\in U$ and $w\in W$. Then
\[
\phi:V(\Gamma)\longrightarrow R\times\{+,-\},\,\,\,\,
u^x\longmapsto(x,+),\quad w^x\longmapsto(x,-)
\]
is a bijection. Put $S=\{h\in R\mid w^h\in\Gamma(u)\}$. Then
$\phi$ is an isomorphism from $\Gamma$ to $\BCay(R,S)$ that identifies
$R$ with $\hat R$. In particular, $1\in S$ if $\{u,w\}$ is an edge.
Thus Theorem~\ref{thm:main} has the following equivalent formulation.

\begin{theorem}\label{thm:main-1}
Let $R$ be a nonabelian simple group and $S\subseteq R$ such that
$|S|$ is prime and $R=\langle S^{-1}S\rangle$. Let
$\Gamma=\BCay(R,S)$ and $A=\Aut(\Gamma)$. If $N_A(\hat R)$ is
edge-transitive on $\Gamma$, then $\hat R$ is a characteristic subgroup
of $A$.
\end{theorem}

Let $\Gamma=\BCay(R,S)$. The normalizer
$N_{\Aut(\Gamma)}(\hat R)$ is described in \cite{Lu03}; see also
\cite{CZFF,ZF}. For $g,h\in R$ and $\alpha\in\Aut(R)$, define
permutations $\kappa_{g,h,\alpha}$ and $\iota_{g,h,\alpha}$ of
$R\times\{+,-\}$ by
\[
\begin{aligned}
(x,+)^{\kappa_{g,h,\alpha}}&=(gx^\alpha,+),&
(x,-)^{\kappa_{g,h,\alpha}}&=(hx^\alpha,-),\\
(x,+)^{\iota_{g,h,\alpha}}&=(gx^\alpha,-),&
(x,-)^{\iota_{g,h,\alpha}}&=(hx^\alpha,+).
\end{aligned}
\]
The following lemma holds.

\begin{lemma}\label{lem:binormalizer}
Let $R$ be a finite group and $S\subseteq R$. Put $\Gamma=\BCay(R,S)$,
$A=\Aut(\Gamma)$ and $N=N_A(\hat R)$. Then
\[
\begin{aligned}
N={}&\{\kappa_{g,h,\alpha}\mid g,h\in R,\ \alpha\in\Aut(R),\ hS^\alpha=Sg\}\\
&\mathbin{\cup}\{\iota_{g,h,\alpha}\mid g,h\in R,\ \alpha\in\Aut(R),\ hS^\alpha=S^{-1}g\},\\
N_{(1,+)}={}&\{\kappa_{1,h,\alpha}\mid h\in R,\ \alpha\in\Aut(R),\ hS^\alpha=S\},\\
N_{(1,-)}={}&\{\kappa_{g,1,\alpha}\mid g\in R,\ \alpha\in\Aut(R),\ S^\alpha=Sg\}.
\end{aligned}
\]
If, in addition, $R=\langle S^{-1}S\rangle$, then $N_{(1,+)}$ and
$N_{(1,-)}$ act faithfully on $\Gamma((1,+))$ and $\Gamma((1,-))$,
respectively.
\end{lemma}

Write
\[
\Aut(R,S)=\{\alpha\in\Aut(R)\mid S^\alpha=S\}.
\]
The holomorph $\mathrm{Hol}(R)=\hat R\Aut(R)$ intersects
$\Aut(\Cay(R,S))$ in $\hat R\Aut(R,S)$; equivalently,
\[
N_{\Aut(\Cay(R,S))}(\hat R)=\hat R\Aut(R,S);
\]
see \cite{Godsil}. The stabilizer of $1$ in this normalizer is
$\Aut(R,S)$. Consequently, $N_{\Aut(\Cay(R,S))}(\hat R)$ is
arc-transitive on $\Cay(R,S)$ if and only if $\Aut(R,S)$ is transitive
on $S$. Combining this observation with Theorem~\ref{thm:main-1}
gives the following corollary.

\begin{corollary}\label{cor:1}
Let $R$ be a nonabelian simple group and $1\notin S\subseteq R$ such
that $R=\langle S^{-1}S\rangle$. Suppose that $|S|$ is prime and
$\Aut(R,S)$ is transitive on $S$. Then
$\hat R\unlhd\Aut(\BCay(R,S))$; in particular,
$\hat R\unlhd\Aut(\Cay(R,S))$.
\end{corollary}
\begin{proof}
Let $\Gamma=\BCay(R,S)$ and $A=\Aut(\Gamma)$. Then
\[
\hat R\Aut(R,S)\leqslant\Aut(\Cay(R,S))\leqslant A^+.
\]
By assumption, $\hat R\Aut(R,S)$ is locally arc-transitive on $\Gamma$.
Since $R=\langle S^{-1}S\rangle$, the graph $\Gamma$ is connected by
\cite[Lemma~2.3(iii)]{DuXu}. Moreover, $\Gamma$ has prime valency and
$\hat R\Aut(R,S)\leqslant N_A(\hat R)$. Theorem~\ref{thm:main-1}
therefore gives $\hat R\unlhd A$, and hence
$\hat R\unlhd\Aut(\Cay(R,S))$.
\end{proof}

The following result is a consequence of Corollary~\ref{cor:1}.

\begin{corollary}\label{cor:2}
Let $R$ be a nonabelian simple group and $1\notin S\subseteq R$ such
that $S^{-1}=S$ and $R=\langle S\rangle$. Suppose that $|S|$ is prime
and $\hat R\Aut(R,S)$ is edge-transitive on $\Cay(R,S)$. Then $|S|$
is odd, $S$ consists of involutions, $R=\langle SS\rangle$, and
$\hat R\unlhd\Aut(\BCay(R,S))$; in particular,
$\hat R\unlhd\Aut(\Cay(R,S))$.
\end{corollary}
\begin{proof}
If $|S|=2$, then the inverse-closed set $S$ either has the form
$\{x,x^{-1}\}$ with $x\ne x^{-1}$ or consists of two distinct
involutions. In the first case $R=\langle S\rangle$ is cyclic, and in
the second case it is dihedral. Both cases contradict nonabelian
simplicity. Thus $|S|$ is an odd prime. Since $S^{-1}=S$ and $1\notin S$,
the set $S$ contains an involution.

Note that a vertex- and edge-transitive group on a graph of odd valency is
arc-transitive. Then $\hat R\Aut(R,S)$ is arc-transitive on
$\Cay(R,S)$, and so $\Aut(R,S)$ is transitive on
$S$. All elements of $S$ therefore have the same order, so every
element of $S$ is an involution.

Put $R_0=\langle SS\rangle$. Then $R_0\ne 1$. Since $S=S^{-1}$, conjugation by every
element of $S$ preserves $R_0$, and hence $R_0\unlhd R$. Now the simplicity of $R$ gives
$R=R_0=\langle SS\rangle=\langle S^{-1}S\rangle$. Then the result   follows from
Corollary~\ref{cor:1}.
\end{proof}

\begin{corollary}\label{cor:2-1}
Under the hypotheses of Corollary~\ref{cor:2}, we have
\[
\Aut(\BCay(R,S))=
\bigl(\hat R\Aut(R,S)\bigr)\times\langle\iota_{1,1,1}\rangle.
\]
\end{corollary}
\begin{proof}
Let $\Gamma=\BCay(R,S)$ and $A=\Aut(\Gamma)$. By
Corollary~\ref{cor:2}, $A=N_A(\hat R)$. Since $\hat R$ is regular on
each part, $A^+=\hat R A_{(1,+)}$. Lemma~\ref{lem:binormalizer} gives
\[
A_{(1,+)}=
\{\kappa_{1,h,\alpha}\mid h\in R,\ \alpha\in\Aut(R),\ hS^\alpha=S\}.
\]
Let $\kappa_{1,h,\alpha}\in A_{(1,+)}$. Since every element of $S$, and hence
of $S^\alpha=h^{-1}S$, is an involution, for every $x\in S$, we have
\[
h=((h^{-1}x)x)^{-1}=x^{-1}(h^{-1}x)^{-1}=xh^{-1}x.
\]
Consequently, every element of $S$ inverts $h$, and every element of
$SS$ centralizes $h$. Since $R=\langle SS\rangle$, we have $h\in Z(R)$.
As $R$ is nonabelian simple, $h=1$. It follows that
$A_{(1,+)}=\Aut(R,S)$ and $A^+=\hat R\Aut(R,S)$. Finally,
$\iota_{1,1,1}\in A$ interchanges the two parts and centralizes $A^+$.
Since $|A:A^+|=2$, the stated direct product decomposition follows.
\end{proof}

\vskip 20pt

\section{Examples of semisymmetric graphs}\label{sec:examples}

A connected regular graph $\Gamma$ is called semisymmetric if $\Aut(\Gamma)$ is
edge-transitive but not vertex-transitive; equivalently, it is locally
arc-transitive but not vertex-transitive on $\Gamma$. Folkman
\cite{Folkman} initiated the study of semisymmetric graphs and posed
eight open problems. These problems stimulated research on their
construction and classification; see, for example,
\cite{Bouwer2,768-ss,Du99,DuXu,IOI,LiuLu,2Odd}. In this section we use
Lemma~\ref{lem:binormalizer} and Corollary~\ref{cor:1} to construct a
class of semisymmetric graphs.

\begin{corollary}\label{cor:semisymmetry}
Let $R$ be a nonabelian simple group and $1\notin S\subseteq R$ such
that $R=\langle S^{-1}S\rangle$. Suppose that $|S|$ is prime and
$\Aut(R,S)$ is transitive on $S$. Then $\BCay(R,S)$ is semisymmetric
if and only if $S^{-1}\ne xS^\alpha$ for all $x\in R$ and
$\alpha\in\Aut(R)$.
\end{corollary}
\begin{proof}
Let $\Gamma=\BCay(R,S)$ and $A=\Aut(\Gamma)$. Then $\Gamma$ is
connected, and $\hat R\Aut(R,S)$ is locally arc-transitive on $\Gamma$.
Hence $\Gamma$ is semisymmetric if and only if $A$ is not
vertex-transitive, that is, $A=A^+$. By Corollary~\ref{cor:1} and
Lemma~\ref{lem:binormalizer},
\[
A\setminus A^+=
\{\iota_{g,h,\alpha}\mid g,h\in R,\ \alpha\in\Aut(R),\ hS^\alpha=S^{-1}g\}.
\]
If $S^{-1}=xS^\alpha$ for some $x\in R$ and $\alpha\in\Aut(R)$, then
$\iota_{1,x,\alpha}\in A\setminus A^+$. Conversely, if
$\iota_{g,h,\alpha}\in A\setminus A^+$, then $hS^\alpha=S^{-1}g$.
Let $\beta$ be the inner automorphism induced by $g^{-1}$ and put
$x=hg^{-1}$. Then $S^{-1}=xS^{\alpha\beta}$, and the result follows.
\end{proof}

We now apply Corollary~\ref{cor:semisymmetry} to the alternating group
$\mathrm A_p$, where $p>5$ is prime. Label the points of its natural
action by the elements of $\mathbb F_p=\mathbb Z/p\mathbb Z$, and let
$\rho\in\mathrm A_p$ be the $p$-cycle defined by
\[
\rho:\mathbb F_p\longrightarrow\mathbb F_p,\,\,\,\, u\longmapsto u+1.
\]
For every integer $i$ and every $u\in\mathbb F_p$, we have
$u^{\rho^i}=u+i$, where $i$ is interpreted modulo $p$. Let $r>2$ be
prime with $p>2r$, and let $\delta\in\mathrm A_p$ be an $r$-cycle. Put
\[
\Delta_r=\{\delta^{\rho^i}\mid 0\leqslant i\leqslant p-1\}.
\]
We denote by $\operatorname{supp}(\sigma)$ the support of
$\sigma\in\mathrm S_p$.

\begin{lemma}\label{N(D_r)}
We have
\[
N_{\mathrm S_p}(\Delta_r)=\langle\rho\rangle{:}\langle\delta^*\rangle,
\]
where $\delta^*$ is either the identity or a product of
$\frac{p-1}{r}$ disjoint $r$-cycles, one of which is $\delta$.
\end{lemma}
\begin{proof}
Put $P=\langle\rho\rangle$ and $N=N_{\mathrm S_p}(\Delta_r)$. Since
$|P|=p$, the set $\Delta_r$ has size $1$ or $p$. In the first case
$\rho$ would centralize $\delta$ and fix $\operatorname{supp}(\delta)$ set-wise, contrary to the transitivity of $P$. Thus $|\Delta_r|=p$,
and clearly $P\leqslant N$.

For distinct $u,v\in\mathbb F_p$, define
\[
m(u,v)=|\{\gamma\in\Delta_r\mid u^\gamma=v\}|.
\]
Then these numbers are invariant under $N$, that is, $m(u^\sigma,v^\sigma)=m(u,v)$ for all $u,v\in\mathbb F_p$ and $\sigma\in N$. 
Note that each $\gamma\in\Delta_r$ is an $r$-cycle and hence contributes exactly
$r$ to $\sum_{u\ne v}m(u,v)$, one for each point in its support.
Suppose that $N$ is doubly transitive. Then all $m(u,v)$
  have a common integral value, say $\lambda$.  Counting the triples $(u,v,\gamma)$ with distinct $u,v\in\mathbb F_p$ and $\gamma\in \Delta_r$, we have 
\[
pr=\sum_{u\ne v}m(u,v)=\lambda p(p-1),
\]
contrary to $0<r<p-1$. Therefore, Burnside's theorem for groups of prime degree
gives
\[
P\unlhd N\leqslant N_{\mathrm S_p}(P)=\mathrm{AGL}_1(p);
\]
see \cite[Theorem~3.5B]{Dixon}.

Since $P$ acts regularly on $\Delta_r$ by conjugation, we have
$N=P{:}C_N(\delta)$. Restriction to $\operatorname{supp}(\delta)$
embeds $C_N(\delta)$ into the cyclic group generated by the restriction
of $\delta$: the kernel fixes at least two points and is trivial
because its elements are affine permutations. Thus $|C_N(\delta)|$
is either $1$ or $r$. In the first case take $\delta^*=1$.

Suppose now that $|C_N(\delta)|=r$.
Choose a generator $c$ of $C_N(\delta)$ and write
\[
c|_{\operatorname{supp}(\delta)}
=\delta^k|_{\operatorname{supp}(\delta)},\,\,\,\, k\not\equiv0\pmod r.
\]
Choose $e$ such that $ke\equiv1\pmod r$ and put $\delta^*=c^e$.
This is a nonidentity affine permutation of order $r$, so it has
exactly one fixed point and all its other orbits have length $r$.
Consequently, it is a product of $(p-1)/r$ disjoint $r$-cycles, one
of which is $\delta$, and $N=P{:}\langle\delta^*\rangle$.
\end{proof}

Let
\[
t=\frac{1}{r}\sum_{x\in\operatorname{supp}(\delta)}x\in\mathbb F_p,
\,\,\,\,
\tau:\mathbb F_p\longrightarrow\mathbb F_p,\quad u\longmapsto-u+2t.
\]
Then $\tau$ is an involution and $\rho^\tau=\rho^{-1}$. Put
\[
\tau_i=\tau\rho^i\quad(0\leqslant i\leqslant p-1),
\,\,\,\,
I_r=\{\sigma\in\mathrm S_p\mid\Delta_r^\sigma=\Delta_r^{-1}\}.
\]

\begin{lemma}\label{I(D_r)}
We have
\[
I_r=
\begin{cases}
\{\tau_i\mid0\leqslant i\leqslant p-1\},&\text{if }\delta^\tau=\delta^{-1},\\
\emptyset,&\text{otherwise}.
\end{cases}
\]
\end{lemma}
\begin{proof}
Suppose that $I_r\ne\emptyset$ and choose $\sigma\in I_r$. Then
\[
N_{\mathrm S_p}(\Delta_r)^\sigma
=N_{\mathrm S_p}(\Delta_r^{-1})=N_{\mathrm S_p}(\Delta_r).
\]
By Lemma~\ref{N(D_r)}, $\langle\rho\rangle$ is the unique Sylow
$p$-subgroup of $N_{\mathrm S_p}(\Delta_r)$, and hence is
characteristic in that group. It follows that
\[
\sigma\in N_{\mathrm S_p}(\langle\rho\rangle)=\mathrm{AGL}_1(p).
\]
Write $\delta^\sigma=(\delta^{-1})^{\rho^i}$ for some $i$, and put
$\beta=\sigma\rho^{-i}$. Then $\delta^\beta=\delta^{-1}$ and
$\beta\ne1$. Write $\delta=(u_0\,u_1\,\ldots\,u_{r-1})$. The equality
\[
(u_0^\beta\,u_1^\beta\,\ldots\,u_{r-1}^\beta)
=\delta^{-1}=(u_{r-1}\,u_{r-2}\,\ldots\,u_1\,u_0)
\]
implies that there is $k\in\{0,1,\ldots,r-1\}$ such that
$u_j^\beta=u_{k-j}$ for $0\leqslant j\leqslant r-1$, with subscripts
taken modulo $r$. Thus $\beta^2$ fixes $\operatorname{supp}(\delta)$
pointwise. Since $\beta^2\in\mathrm{AGL}_1(p)$ fixes at least two
points, $\beta^2=1$. As $p$ is odd and $\beta\ne1$, there is
$v\in\mathbb F_p$ such that $u^\beta=-u+v$ for all $u\in\mathbb F_p$.
For $0\leqslant j\leqslant r-1$, we have $v=u_j+u_{k-j}$, and hence
\[
v=\frac{1}{r}\sum_{j=0}^{r-1}(u_j+u_{k-j})=2t.
\]
Therefore $\beta=\tau$, $\sigma=\tau\rho^i=\tau_i$ and
$\delta^\tau=\delta^{-1}$.

Conversely, suppose that $\delta^\tau=\delta^{-1}$. Since
$\rho^\tau=\rho^{-1}$, for every $0\leqslant i\leqslant p-1$ we have
\[
\begin{aligned}
\Delta_r^{\tau_i}
&=\{(\delta^{\rho^j})^{\tau\rho^i}\mid0\leqslant j\leqslant p-1\}\\
&=\{(\delta^{-1})^{\rho^{i-j}}\mid0\leqslant j\leqslant p-1\}
=\Delta_r^{-1}.
\end{aligned}
\]
Together with the first part, this proves the lemma.
\end{proof}

Suppose now that $r=3$, so that $\delta=(u_0\,u_1\,u_2)$ is a
$3$-cycle. Then $\delta^\tau=\delta^{-1}$ if and only if $\tau$
fixes one point of $\operatorname{supp}(\delta)$, say $u_1$, and
interchanges the other two. In this case, $u_0+u_2=2u_1$. We therefore
obtain the following consequence of Lemma~\ref{I(D_r)}.

\begin{corollary}\label{I_3}
Let $p>5$ be prime. Then $I_3\ne\emptyset$ if and only if there exist
$u,d\in\mathbb F_p$ with $d\ne0$ such that
$\operatorname{supp}(\delta)=\{u-d,u,u+d\}$.
\end{corollary}

\begin{lemma}\label{A_p}
Let $p>5$ be prime. Then $\langle\Delta_3^{-1}\Delta_3\rangle=\mathrm A_p$.
If $\sigma\in\mathrm A_p$ is such that $\sigma\Delta_3$ consists of
$3$-cycles, then $\sigma=1$.
\end{lemma}
\begin{proof}
Put $D=\operatorname{supp}(\delta)$ and $\delta_i=\delta^{\rho^i}$.
Then $\operatorname{supp}(\delta_i)=D+i\subseteq\mathbb F_p$, where
$i$ is interpreted modulo $p$. For $i\ne0$, put
$\ell_i=|D\cap(D+i)|$. Since $D$ is a nonempty proper subset of
$\mathbb F_p$, it cannot be invariant under a nontrivial translation;
thus $\ell_i<3$. The map $u\mapsto u-i$ is a bijection from
$D\cap(D+i)$ to $D\cap(D-i)$, so $\ell_i=\ell_{-i}$. Also,
\[
D\cap(D+i)\longrightarrow\{(u,v)\in D\times D\mid u-v=i\},
\,\,\,\, u\longmapsto(u,u-i),
\]
is a bijection. Counting ordered pairs of distinct points of $D$ gives
\[
6=\sum_{i\ne0}\ell_i
=\sum_{i=1}^{(p-1)/2}(\ell_i+\ell_{-i})
=2\sum_{i=1}^{(p-1)/2}\ell_i.
\]
Since $\ell_i\in\{0,1,2\}$ and $6$ is not divisible by $4$, we have
$\ell_i=1$ for some $i\ne0$. Hence two elements $\delta_j,\delta_k$
of $\Delta_3$ have supports meeting in exactly one point, and
$\delta_j^{-1}\delta_k\in\Delta_3^{-1}\Delta_3$ is a $5$-cycle.

Put $H=\langle\Delta_3^{-1}\Delta_3\rangle$. Since $\langle \rho\rangle$ normalizes
$H$ and is transitive on $\mathbb F_p$, it permutes the $H$-orbits
transitively. These orbits therefore have equal size. As $H\ne1$ and
$p$ is prime, $H$ is transitive and hence primitive. If $p>7$, Jordan's
theorem gives $H=\mathrm A_p$; see
\cite[p.~84, Theorem~3.3E]{Dixon}.

It remains to consider $p=7$. Put $m=|H|$. Then $35\mid m$ and
$m\mid2520$. Let $P$ and $Q$ be Sylow $7$- and Sylow $5$-subgroups
of $H$, respectively. We have $|N_{\mathrm A_7}(P)|=21$ and
$|N_{\mathrm A_7}(Q)|=20$. Thus $|N_H(P)|\in\{7,21\}$ and
$|N_H(Q)|\mid20$. The number of elements of order $5$ in $H$ is at
least $4m/20$. If $|N_H(P)|=7$, the number of elements of order $7$
is $6m/7$, giving
\[
m>\frac{4m}{20}+\frac{6m}{7}=\frac{37m}{35}>m,
\]
a contradiction. Hence $|N_H(P)|=21$. The number $n_7$ of Sylow
$7$-subgroups of $H$ satisfies
\[
n_7=\frac{m}{21},\,\,\,\, n_7\mid120,\,\,\,\, 5\mid n_7,
\,\,\,\, n_7\equiv1\pmod7.
\]
It follows that $n_7=15$ or $120$, and hence $m=315$ or $2520$.
If $m=315$, then $|N_H(Q)|=5$, so $H$ has exactly $63$ Sylow
$5$-subgroups, contrary to $63\not\equiv1\pmod5$. Therefore
$m=2520$ and $H=\mathrm A_7$. This proves the first assertion.

For the second assertion, suppose that $\sigma\delta_i$ is a
$3$-cycle for every $i$. Recall that ${\rm supp}(\delta_i)=D+i$, where $0\leqslant i\leqslant p-1$, and so    $u\in \mathbb{F}_p$ lies in some $D+i$ if and only if $u-i\in D$. Thus each element of $\mathbb F_p$ lies in exactly three of the $p$ translates $D+i$ of $D$. Then there are exactly $3| {\rm supp}(\sigma)|$  ordered pairs $(u,{\rm supp}(\delta_i))$ with $u\in {\rm supp}(\sigma)\cap {\rm supp}(\delta_i)$. We have
\[3| {\rm supp}(\sigma)|=\sum_i|{\rm supp}(\sigma)\cap{\rm supp}(\delta_i)|.\]
Note that   ${\rm supp}(\sigma\delta_i)$ contains
 the symmetric difference ${\rm supp}(\sigma)\mathbin{\triangle}{\rm supp}(\delta_i)$ of  ${\rm supp}(\sigma)$ and  ${\rm supp}(\delta_i)$.
 Consequently,
\[
3\geqslant |{\rm supp}(\sigma)\mathbin{\triangle}{\rm supp}(\delta_i)|
=|{\rm supp}(\sigma)|+3-2|{\rm supp}(\sigma)\cap {\rm supp}(\delta_i)|,
\]
 so $|{\rm supp}(\sigma)|\leqslant2|{\rm supp}(\sigma)\cap {\rm supp}(\delta_i)|\leqslant 6$.   Then
 \[p|{\rm supp}(\sigma)|\leqslant \sum_i2|{\rm supp}(\sigma)\cap {\rm supp}(\delta_i)|=6|{\rm supp}(\sigma)|,\]
  so $p|{\rm supp}(\sigma)|\leqslant 6|{\rm supp}(\sigma)|$.
Since $p>5$ is prime, 
we have $p\geqslant7$, and hence $|{\rm supp}(\sigma)|=0$. Thus $\sigma=1$, and the lemma follows.
\end{proof}

\begin{theorem}\label{th-A_p}
Let $p>5$ be prime, and let $\Gamma=\BCay(\mathrm A_p,\Delta_3)$. Then
$\Gamma$ is semisymmetric if and only if $\operatorname{supp}(\delta)$
is not a three-term arithmetic progression in $\mathbb F_p$.
\end{theorem}
\begin{proof}
Put $A=\Aut(\Gamma)$. The group $\Aut(\mathrm A_p,\Delta_3)$
contains the inner automorphism induced by $\rho$, and hence is
transitive on $\Delta_3$. It follows that $A$ is locally
arc-transitive on $\Gamma$. By Lemma~\ref{A_p},
$\langle\Delta_3^{-1}\Delta_3\rangle=\mathrm A_p$, so $\Gamma$ is
connected. Thus $\Gamma$ is semisymmetric if and only if $A$ is not
vertex-transitive on $\Gamma$.

Since $p>5$, we identify $\Aut(\mathrm A_p)$ with $\mathrm S_p$
acting by conjugation. Suppose that $A$ is vertex-transitive on
$\Gamma$. By Corollary~\ref{cor:semisymmetry}, there exist
$\sigma\in\mathrm A_p$ and $\alpha\in\mathrm S_p$ such that
$\Delta_3^{-1}=\sigma\Delta_3^\alpha$. Therefore
$\sigma^{-1}\Delta_3^{-1}=\Delta_3^\alpha$ consists of $3$-cycles.
Since
\[
\Delta_3^{-1}
=\{(\delta^{-1})^{\rho^i}\mid0\leqslant i\leqslant p-1\},
\]
applying Lemma~\ref{A_p} to $\Delta_3^{-1}$ and $\sigma^{-1}$ gives
$\sigma^{-1}=1$. Hence $\Delta_3^\alpha=\Delta_3^{-1}$.

Conversely, if $\Delta_3^\alpha=\Delta_3^{-1}$ for some
$\alpha\in\Aut(\mathrm A_p)$, then
$\iota_{1,1,\alpha}\in A$ interchanges the two parts of $\Gamma$.
Since $\hat{\mathrm A}_p$ is transitive on each part, $A$ is
vertex-transitive. Consequently, $\Gamma$ is semisymmetric if and
only if $I_3=\emptyset$. Then the theorem follows from Corollary~\ref{I_3}.
\end{proof}

\end{document}